\documentclass[a4paper,12pt,reqno]{article}
\usepackage{graphicx}
\usepackage{latexsym}
\usepackage{amssymb}
\usepackage{amsmath}
\usepackage{layout}  
\usepackage{graphicx}
\usepackage{amssymb}
\usepackage{amssymb}
\usepackage[T1]{fontenc}
\usepackage{latexsym}
\usepackage{xypic}
\usepackage{eufrak}
\usepackage{euscript}
\usepackage{amsfonts,amsmath}
\usepackage{verbatim}
\usepackage{fancyhdr}
\usepackage[english]{babel}
\usepackage{mathrsfs}

\newtheorem{thm}{Theorem}
\newtheorem{defi}{Definition}
\newtheorem{lemma}{Lemma}
\newtheorem{pro}{Proposition}
\newtheorem{cor}{Corollary}

\def\cal{\mathcal}
\def\scr{\mathscr}

\renewcommand{\geq}{\geqslant}
\def\leq{\leqslant}

\newcommand{\N}{\mathbb{N}}

\newcommand{\R}{\mathbb{R}}
\newcommand{\Q}{{\rm \bf Q}}
\newcommand{\bP}{{\rm \bf P}}

\newcommand{\fr}{\EuFrak}

\def\qed{ \hfill \vrule width.25cm height.25cm depth0cm\smallskip}

\def\esp{{\rm E}}
\def\var{{\mathbb{Var}}}

\def\e{\varepsilon}

\def\cal{\mathcal}
\def\1{{\mathbf{1}}}

\def\sg{\sigma}

\def\oT{{[0{,}T]}}

\def\ffi{{\varphi}}

\def\1{{\mathbf{1}}}
\def\0.5{{\frac{1}{2}}}
\def\var{{\rm{Var}}}
\def\cov{{\rm{Cov}}}

\def\sg{\sigma}

\newcommand{\ignore}[1]{}
\begin{document}

\title{Differentiating $\sg$-fields for Gaussian and shifted Gaussian processes}
\author{
S\'ebastien Darses\footnote{
LPMA, Universit{\'e} Paris 6,
Bo{\^\i}te courrier 188, 4 Place Jussieu, 75252 Paris Cedex 05, France,
{\tt sedarses@ccr.jussieu.fr}
},\,\, 
Ivan Nourdin\footnote{
LPMA, Universit{\'e} Paris 6,
Bo{\^\i}te courrier 188, 4 Place Jussieu, 75252 Paris Cedex 05, France,
{\tt nourdin@ccr.jussieu.fr}
}\,\,
and\,\, Giovanni Peccati\footnote{
LSTA, Universit{\'e} Paris 6,
Bo{\^\i}te courrier 188, 4 Place Jussieu, 75252 Paris Cedex 05, France,
{\tt giovanni.peccati@gmail.com}
}
}
\maketitle

\begin{abstract}
We study the notions of differentiating and non-differentiating
$\sg$-fields in the general framework of (possibly drifted)
Gaussian processes, and characterize their invariance properties
under equivalent changes of probability measure. As an
application, we investigate the class of stochastic derivatives
associated with shifted fractional Brownian motions. We finally
establish conditions for the existence of a jointly measurable
version of the differentiated process, and we outline a general
framework for stochastic embedded equations.
\end{abstract}

\vspace{1cm}

\section{Introduction}
\indent

{\normalsize \label{intro} }

Let $X$ be the solution of the stochastic differential equation
$X_t=X_0+\int_0^t\sigma(X_s)dB_s + \int_0^tb(X_s)ds,\, t\in[0,T]$,
where $\sigma,b:\R\rightarrow\R$ are suitably regular functions and
$B$ is a standard Brownian motion, and denote by $\scr{P}^X_t$ the
$\sigma$-field generated by \{$X_s$, $s\in[0,t]$\}. Then, the
following quantity:
\begin{equation}\label{gene}
h^{-1}{\rm E}\left[f(X_{t+h}) - f(X_t)|\scr{P}^X_t\right]
\end{equation}
converges (in probability and for $h\downarrow 0$) for every smooth
and bounded function $f$. This existence result is the key to define
one of the central operators in the theory of diffusion processes:
the \textit{infinitesimal generator} $L$ of $X$, which is given by
$Lf(x)=b(x)\frac{df}{dx}(x)+\frac{1}{2}
\sigma(x)^2\frac{d^2f}{dx^2}(x)$ (the domain of $L$ contains all
regular functions $f$ as above). Note that the limit in (\ref{gene})
is taken conditionally to the past of $X$ before $t$; however, due
to the Markov property of $X$, one may as well replace $\scr{P}^X_t$
with the $\sigma$-field $\sg\{X_t\}$ generated by $X_t$. On the
other hand, under rather mild conditions on $b$ and $\sigma$, one
can take $f={\rm Id}$ in (\ref{gene}), so that the limit still
exists and coincides with the natural definition of the \textit{mean
velocity} of $X$ at $t$ (the reader is referred to Nelson's
\textsl{dynamical theory of Brownian diffusions}, as developed
\textit{e.g.} in \cite{nels}, for more results in this direction --
see also \cite{carl} for a recent survey).

\smallskip

In this paper we are concerned with the following question:
\textsl{is it possible to obtain the existence, and to study the
nature, of limits analogous to (\ref{gene}), when $X$ is neither a
Markov process nor a semimartingale?} We will mainly focus on the
case where $X$ is a (possibly shifted) Gaussian random process and
$f={\rm Id}$ (the case of a non-linear and smooth $f$ will be
investigated elsewhere). The subtleties of the problem are better
appreciated through an example. Consider for instance a fractional
Brownian motion (fBm) $B$ of Hurst index $H \in (1/2,1)$, and
recall that $B$ is neither Markovian nor a semimartingale (see
\textit{e.g.} \cite{n}). Then, the quantity $h^{-1} {\rm
E}[B_{t+h}-B_t|B_t]$ converges in $L^2(\Omega)$ (as
$h\downarrow0$), while the quantity $h^{-1} {\rm
E}[B_{t+h}-B_t|\mathscr{P}^B_t]$ does not admit a limit in
probability. More to the point, similar properties can be shown to
hold also for suitably regular solutions of stochastic
differential equations driven by $B$ (see \cite{dn} for precise
statements and proofs).

\smallskip

To address the problem evoked above, we shall mainly use the
notion of {\it differentiating $\sg$-field} introduced in
\cite{dn}: if $Z$ is a process defined on a probability space
$(\Omega,\scr F,\bP)$, we say that a $\sg$-field $\scr
G\subset\scr F$ is \textit{differentiating} for $Z$ at $t$ if
\begin{equation}\label{tribus}
h^{-1}{\rm E}\left[Z_{t+h}-Z_t|\scr G\right]
\end{equation}
converges in some topology, when $h$ tends to $0$. When it exists,
the limit of (\ref{tribus}) is noted $D^{\scr G}Z_t$, and it is
called the \textit{stochastic derivative} of $Z$ at $t$ with
respect to $\scr G$. Note that if a sub-$\sigma$-field
$\mathscr{G}$ of $\mathscr{F}$ is not differentiating, one can
implement two \textquotedblleft strategies\textquotedblright\ to
make (\ref{tribus}) converge: either one replaces $\mathscr{G}$
with a differentiating sub-$\sigma$-field $\mathscr{H}$, or one
replaces $h^{-1}$ with $h^{-\alpha}$ with $0<\alpha<1$. In
particular, the second strategy pays dividends when a
non-differentiating $\sigma$-field $\scr G$ is \textsl{too poor},
in the sense that $\scr G$  does not contain sufficiently good
differentiating $\sg$-fields. We will see that this is exactly the
case for a fBm $B$ with index $H<1/2$, when $\scr G$ is generated
by $B_s$ for some $s>0$.

\smallskip

The aim of this paper is to give a precise characterization of the
classes of differentiating and non differentiating $\sg$-fields
for Gaussian and shifted Gaussian processes. We will
systematically investigate their mutual relations, and pay special
attention to their invariance properties under equivalent changes
of probability measure.

\smallskip

The paper is organized as follows. In Sections 2 and 3 we
introduce several notions related to the concept of
differentiating $\sg$-field, and give a characterization of
differentiating and non differentiating $\sg$-fields in a Gaussian
framework. In Section 4 we prove some invariance properties of
differentiating $\sg$-fields under equivalent changes of
probability measure. Notably, we will be able to write an explicit
relation between the stochastic derivatives associated with
different probabilities. We will illustrate our results by
considering the example of shifted fractional Brownian motions,
and we shall pinpoint different behaviors when the Hurst index is,
respectively, in $(0,1/2)$ and in $(1/2,1)$ . In Section 5 we
establish fairly general conditions, ensuring the existence of a
jointly measurable version of the differentiated process induced
by a collection of differentiating $\sg$-fields. Finally, in
Section 6 we outline a general framework for \textit{embedded
ordinary stochastic differential equations} (as defined in
\cite{cd}) and we analyze a simple example.

\section{\protect\normalsize Preliminaries on stochastic derivatives}
{\normalsize \label{sec:1}}

Let $(Z_t)_{t\in [0,T]}$ be a stochastic process defined on a
probability space $(\Omega,\mathscr{F}, \bP)$. In the sequel, we
will
always assume that $Z_t\in \mathrm{L}^2(\Omega,%
\mathscr{F}, \bP)$ for every $t\in [0,T]$. It will also be
implicit that each $\sigma$-field we consider is a
sub-$\sigma$-field of $\mathscr{F}$; analogously, given a
$\sigma$-field $\mathscr{H}$, the notation $\scr G\subset\scr H$
will mean that $\scr G$ is a sub-$\sigma$-field of $\mathscr{H}$.
For every $t\in(0,T)$ and every $h\neq 0$ such that $t+
h\in(0,T)$, we set
\begin{equation*}
\Delta_hZ_t=\frac{Z_{t+h}-Z_t}{h}.
\end{equation*}
For the rest of the paper, we will use the letter $\tau$ as a
generic symbol to indicate a topology on the class of real-valued
and $\mathscr{F}$-measurable random variables. For instance,
$\tau$ can be the topology induced either by the a.s. convergence,
or by the ${\rm L}^p$ convergence ($p\geq1$), or by the
convergence in probability, or by both a.s. and ${\rm L}^p$
convergences, in which cases we shall write, respectively,
$$\tau={\rm a.s.},\quad\tau= {\rm L}^p,\quad\tau={\rm
proba},\quad\tau={\rm L}^p\star {\rm a.s.}.$$ Note that, when no
further specification is provided, any convergence is tacitly
defined with respect to the reference probability measure $\bP$.

\begin{defi}\label{defi1}  Fix $t\in(0,T)$ and let $\mathscr{G}\subset\scr{F}$. We say that $\mathscr{G}$
{\rm $\tau$-differentiates} $Z$ at $t$ if
\begin{equation}\label{def3}
\mbox{ $\mathrm{E}[\Delta_hZ_t\,|\,\mathscr{G}]$ converges w.r.t.
$\tau$ when  $h\to 0$. }
\end{equation}
In this case, we  define the so-called {\rm $\tau$-stochastic
derivative} of $Z$ w.r.t. $\mathscr{G}$ at $t$ by
\begin{equation}\label{der}
D^{\mathscr{G}}_\tau Z_t = \tau\mbox{-}\lim_{h\to 0}\mathrm{E}[
\Delta_hZ_t\,|\,\mathscr{G}].
\end{equation}
If the limit in (\ref{def3}) does not exist, we say that
$\mathscr{G}$ {\rm does not $\tau$-differentiate} $Z$ at $t$. If
there is no risk of ambiguity on the topology $\tau$, we will
write $D^{\mathscr{G}} Z_t$ instead of $D^{\mathscr{G}}_\tau Z_t$
to simplify the notation.
\end{defi}

\textbf{Remark.} When $\tau={\rm a.s.}$ (i.e., when $\tau$ is the
topology induced by a.s. convergence), equation (\ref{def3}) must
be understood in the following sense (note that, in (\ref{def3}),
$t$ acts as a fixed parameter): there exists a jointly measurable
application $(\omega,h)\mapsto q(\omega,h)$, from $\Omega \times
(-\e,\e)$ to $\R$, such that (i) $q(\cdot,h)$ is a version of
$\mathrm{E}[\Delta_hZ_t\,|\,\mathscr{G}]$ for every fixed $h$, and
(ii) there exists a set $\Omega' \subset \Omega$, of
$\bP$-probability one, such that $q(\omega,h)$ converges, as $h\to
0$, for every $\omega \in \Omega'$. An analogous remark applies to
the case $\tau={\rm L}^p\star {\rm a.s.}$ ($p\geq1$).

\smallskip

{\normalsize
The set of all $\sigma$-fields that $\tau$-differentiate $Z$ at time $%
t$ is denoted by $\mathscr{M}^{(t),\tau}_{Z}$. Intuitively, one
can say that the more $\mathscr{M}^{(t),\tau}_Z$ is large, the
more $Z$ is regular at time $t$.
For instance, one has clearly that $\{\emptyset,\Omega\} \in \mathscr{M}%
^{(t),\tau}_{Z}$ if, and only if, the application $s\mapsto
\mathrm{E}(Z_s)$ is differentiable at
time $t$. On the other hand, $\mathscr{F} \in \mathscr{M}%
^{(t),\tau}_{Z}$ if, and only if, the random function $s\mapsto Z_s$ is $\tau$-differentiable at time $t$%
.}

\smallskip

Before introducing some further definitions, we shall illustrate the
above notions by a simple example involving the ${\rm L}^2\star{\rm
a.s.}$-topology. Assume that $Z=(Z_t)_{t\in[0,T]}$ is a Gaussian
process such that ${\rm Var}(Z_t)\neq 0$ for every $t\in(0,T]$. Fix
$t\in(0,T)$ and take $\mathscr{G}$ to be the present of $Z$ at a
fixed time $s\in(0,T]$, that is, $\mathscr{G}=\sigma\{Z_s\}$ is the
$\sigma$-field generated by $Z_s$. Since one has, by linear
regression,
$$
{\rm E}[\Delta_hZ_t\,|\,\mathscr{G}]=\frac{{\rm
Cov}(\Delta_hZ_t,Z_s)}{{\rm Var}(Z_s)}\,Z_s,
$$
we immediately deduce that $\mathscr{G}$ differentiates $Z$ at $t$
if, and only if, $$\frac{d}{du}\,\cov(Z_u,Z_s)|_{u=t}$$ exists
(see also Lemma \ref{derivee=cov}). Now, let $\mathscr{H}$ be a
$\sigma$-field such that $\mathscr{H}\subset\mathscr{G}$. Owing to
the projection principle, one can write:
$$
{\rm E}[\Delta_hZ_t|\mathscr{H}]=\frac{{\rm
Cov}(\Delta_hZ_t,Z_s)}{{\rm Var}(Z_s)}\,{\rm E}[Z_s|\mathscr{H}],
$$
and we conclude that
\begin{enumerate}
\item[(A)] If $\mathscr{G}$ differentiates $Z$ at $t$, then it is also the case for any $\mathscr{H}\subset\mathscr{G}$.
\item[(B)] If $\mathscr{G}$ does not differentiates $Z$ at $t$, then any $\mathscr{H}\subset\mathscr{G}$ either
does not differentiates $Z$ at $t$, or (when ${\rm
E}[Z_s|\mathscr{H}]=0$) differentiates $Z$ at $t$ with
$D^{\mathscr{H}} Z_t=0$

\end{enumerate}
The phenomenon appearing in (A) is quite natural, not only in a
Gaussian setting, and it is due to the well-known properties of
conditional expectations: see Proposition \ref{sous-bonne} below.
On the other hand, (B) seems strongly linked to the Gaussian
assumptions we made on $Z$. We shall use fine arguments to
generalize (B) to a non-Gaussian framework, see Sections
\ref{sec:2} and \ref{sec:3} below.

This example naturally leads to the subsequent definitions.

\begin{defi}
{\normalsize Fix $t\in(0,T)$ and  let $\mathscr{G}\subset\scr{F}$.
If $\mathscr{G}$ $\tau$-differentiates $Z$ at $t$ and if we have
$D^{\mathscr{G}}_\tau Z_t=c$ ${\rm a.s.}$ for a certain real
$c\in\R$, we say that $\mathscr{G}$ {\rm $\tau$-degenerates} $Z$ at
$t$. We say that a random variable $Y$ {\rm $\tau$-degenerates} $Z$
at $t$ if the $\sigma$-field $\sigma\{Y\}$ generated by $Y$
$\tau$-degenerates $Z$ at $t$. }
\end{defi}

{\normalsize If $D^{\mathscr{G}}_\tau Z_t\in{\rm L}^2$ (for
instance when we choose $\tau={\rm L}^2$, or $\tau={\rm
L}^2\star{\rm a.s.}$, etc.), the condition on
$D^{\mathscr{G}}_\tau Z_t$ in the previous definition is obviously
equivalent to ${\rm Var}(D^{\mathscr{G}}_\tau Z_t)=0$.
For instance, if $Z$ is a process such that $s\rightarrow \mathrm{E}%
(Z_s)$ is differentiable in $t\in(0,T)$ then $\{\emptyset,\Omega\}$
degenerates $Z$ at $t$. }

\begin{defi}\label{realy}
\label{rnd} {\normalsize Let $t\in(0,T)$ and
$\mathscr{G}\subset\mathscr{F}$. We say that $\mathscr{G}$ {\rm
really does not $\tau$-differentiate} $Z$ at $t$ if $\mathscr{G}$
does not $\tau$-differentiate $Z$ at $t$ and if any
$\mathscr{H}\subset\mathscr{G}$ either $\tau$-degenerates $Z$ at
$t$, or does not $\tau$-differentiate $Z$ at $t$. }
\end{defi}

{\normalsize Consider e.g. the phenomenon described at point (B)
above: the $\sigma$-field $\mathscr{G}\triangleq\sigma\{Z_s\}$
really does not differentiate the Gaussian process $Z$ at $t$
whenever $\frac{d}{du}\,\cov(Z_u,Z_s)|_{u=t}$ does not exist, since
every $\mathscr{H}\subset\mathscr{G}$ either does not differentiate
or degenerates $Z$ at $t$. It is for instance the case when $Z=B$ is
a fractional Brownian motion with Hurst index $H<1/2$ and $s=t$, see
Corollary \ref{cor_thm2}. Another interesting example is given by
the process $Z_t=f_1(t)N_1+f_2(t)N_2$, where
$f_1,f_2:[0,T]\rightarrow\R$ are two deterministic functions and
$N_1,N_2$ are two centered and independent random variables. Assume
that $f_1$ is differentiable at $t\in(0,T)$ but that $f_2$ is not.
This yields that $\mathscr{G}\triangleq \sigma\{N_1,N_2\}$ does not
differentiates $Z$ at $t$. Moreover, one can easily show that
$\mathscr{H}\triangleq\sigma\{N_1\}\subset\mathscr{G}$
differentiates $Z$ at $t$ with $D^{\scr{H}}Z_t=f'_1(t)\,N_1$, which
is not constant in general. Then, although $\mathscr{G}$ does not
differentiate $Z$ at $t$, it does not meet the requirements of
Definition \ref{realy}.

\section{\protect\normalsize Stochastic derivatives and Gaussian processes}

{\normalsize \label{sec:2} }

{\normalsize In this section we mainly focus on Gaussian
processes, and we shall systematically work with the ${\rm L}^2$-
or the ${\rm L}^2\star{\rm a.s.}$-topology, which are quite
natural in this framework. In the sequel we will also omit the
symbol $\tau$ in (\ref{der}), as we will always indicate the
topology we are working with.

Our aim is to establish several relationships between
differentiating and (really) non differentiating $\sigma$-fields
under Gaussian-type assumptions. However, our first result
pinpoints a general simple fact, which also holds in a
non-Gaussian framework, that is: \textsl{any sub-$\sigma$-field of
a differentiating $\sigma$-field is also differentiating}.
\begin{pro}\label{sous-bonne}
Let $Z$ be a stochastic process (not necessarily Gaussian)
such that $Z_t\in{\rm L}^2(\Omega,%
\mathscr{F}, \bP)$ for every $t\in (0,T)$. Let $t\in (0,T)$ be
fixed, and let $\mathscr{G}\subset \scr F$. If $\mathscr{G}$ ${\rm
L}^2$-differentiates $Z$ at $t$, then any
$\mathscr{H} \subset\mathscr{G}$ also ${\rm L}^2$-differentiates $Z$ at $t$. Moreover, we have
\begin{equation}\label{proj}
D^{\mathscr{H}}Z_t ={\rm E}[D^{\mathscr{G}}Z_t|\mathscr{H}].
\end{equation}
\end{pro}

\begin{proof}
We can write, by the projection principle and Jensen inequality:
\begin{eqnarray*}
{\rm E}\left[\big({\rm E}(\Delta_h Z_t\,|\,\mathscr{H}) -{\rm
E}(D^{\mathscr{G}}Z_t\,|\,\mathscr{H})\big)^2\right]
 &=&{\rm E}\left[{\rm E}[{\rm E}(\Delta_h Z_t-D^{\mathscr{G}}Z_t\,|\,\mathscr{G})\,|
\,\mathscr{H}]^2\right]\\
 & \le &{\rm E}\left[\big({\rm E}(\Delta_h Z_t\,|\,\mathscr{G})-D^{\mathscr{G}}Z_t\big)^2\right].
\end{eqnarray*}
So, the ${\rm L}^2$-convergence of ${\rm E}(\Delta_h
Z_t\,|\,\mathscr{H})$ to ${\rm
E}(D^{\mathscr{G}}Z_t\,|\,\mathscr{H})$ as $h\rightarrow 0$ is
obvious.\qed
\end{proof}
} {\normalsize On the other hand, a non differentiating
$\sigma$-field may contain a differentiating $\sigma$-field (for
instance, when the non differentiating $\sigma$-field is generated
both by differentiating and non differentiating random variables).

\smallskip

We now provide a characterization of the really
non-differentiating $\sigma$-fields that are generated by some
subspace of the first Wiener chaos associated with a centered
Gaussian process $Z$, noted $\mathcal{H}_{1}(Z)$. We recall that
$\mathcal{H}_{1}(Z)$ is the ${\rm L}^2$-closed linear vector space
generated by random variables of the type $Z_t$, $t \in[0,T]$.}

\begin{thm}\label{thm1}
Let $I=\{1,2,\ldots,N\}$, with $N\in\N^*\cup\{+\infty\}$ and let
$Z=(Z_t)_{t\in\oT}$ be a centered Gaussian process. Fix $t\in(0,T)$,
and consider a subset $\{Y_i\}_{i\in I}$ of $\mathcal{H}_1(Z)$ such
that, for any $n\in I$, the covariance matrix $M_n$ of
$\{Y_i\}_{1\le i\le n}$ is invertible. Finally, note
$\mathscr{Y}=\sigma\{Y_i, i\in I\}$. Then:
\begin{enumerate}
\item If $\scr{Y}$ ${\rm L}^2$-differentiates $Z$ at $t$, then,
for any $i\in I$, $Y_i$ ${\rm L}^2$-differentiates $Z$ at $t$. If
$N<+\infty$, the converse also holds.
\item Suppose $N<+\infty$. Then $\mathscr{Y}$ really does
not ${\rm L}^2\star {\rm a.s.}$-differentiate $Z$ at $t$ if, and
only if, any finite linear combination of the $Y_i$'s either ${\rm
L}^2\star {\rm a.s.}$-degenerates or does not ${\rm L}^2\star {\rm
a.s.}$-differentiate $Z$ at $t$.
\item Suppose that $N=+\infty$ and that the sequence $\{Y_i\}_{i\in
I}$ is i.i.d.. Write moreover $\mathbf{R}(\mathscr{Y})$ to indicate
the class of all the sub-$\sigma$-fields of $\mathscr{Y}$ that are
generated by rectangles of the type $A_1 \times \cdot\cdot\cdot
\times A_d$, with $A_i\in\sigma\{Y_i\}$, and $d\geq 1$. Then, the
previous characterization holds in a weak sense: if $\mathscr{Y}$
really does not ${\rm L}^2\star {\rm a.s.}$-differentiate $Z$ at
$t$, then every finite linear combination of the $Y_i$'s either
${\rm L}^2\star {\rm a.s.}$-degenerates or does not ${\rm L}^2\star
{\rm a.s.}$-differentiate $Z$ at $t$; on the other hand, if every
finite linear combination of the $Y_i$'s either ${\rm L}^2\star {\rm
a.s.}$-degenerates or does not ${\rm L}^2\star {\rm
a.s.}$-differentiate $Z$ at $t$, then any
$\scr{G}\in\mathbf{R}(\mathscr{Y})$ either ${\rm L}^2\star {\rm
a.s.}$-degenerates or does not ${\rm L}^2\star {\rm
a.s.}$-differentiate $Z$ at $t$.

\end{enumerate}
\end{thm}

The class $\mathbf{R}(\mathscr{Y})$ contains for instance the
$\sg$-fields of the type
$$\mathscr{G}=\sg\{f_1(Y_1),...,f_d(Y_d)\},$$where $d\geq 1$. When
$N=1$, the second point of Theorem \ref{thm1} can be reformulated
as follows (see also the examples discussed in Section \ref{sec:1}
above).
\begin{cor}\label{cor1}
Let $Z=(Z_t)_{t\in\oT}$ be a centered Gaussian process and let
$\mathcal{H}_{1}(Z)$ be its first Wiener chaos. Fix $t\in(0,T)$, as
well as $Y\in\mathcal{H}_1(Z)$, and set $\mathscr{Y}=\sigma\{Y\}$.
Then, $\mathscr{Y}$ does not ${\rm L}^2$-differentiate $Z$ at $t$
(resp. ${\rm L}^2\star {\rm a.s.}$) if, and only if, $\mathscr{Y}$
{\it really} does not ${\rm L}^2$-differentiate $Z$ at $t$ (resp.
${\rm L}^2\star {\rm a.s.}$).
\end{cor}
In particular, when $Z=B$ is a fractional Brownian motion with Hurst
index $H\in(0,1/2)\cup(1/2,1)$, $t$ is a fixed time in $(0,T)$ and
$\mathscr{Y}=\sigma\{B_t\}$ is the present of $B$ at time $t$, we
observe two distinct behaviors, according to the different values of
$H$:
\begin{enumerate}
\item[(a)] If $H>1/2$, then $\mathscr{Y}$ ${\rm L}^2\star{\rm a.s.}$-differentiates $B$ at $t$ and
it is also the case for any $\mathscr{Y}_0\subset\mathscr{Y}$.
\item[(b)] If $H<1/2$, then $\mathscr{Y}$ {\it really} does not ${\rm L}^2\star{\rm a.s.}$-differentiate $B$ at $t$.
\end{enumerate}
Indeed, (a) and (b) are direct consequences of Proposition
\ref{sous-bonne}, Corollary \ref{cor1} and the equality
$$
{\rm E}[\Delta_hB_t|B_t]=\frac{(t+h)^{2H}-t^{2H}-|h|^{2H}}{2\,t^{2H}h}\,B_t,
$$
which is immediately verified by a Gaussian linear regression.

Note that \cite[Theorem 22]{dn} generalizes (a) to the case of
fractional diffusions. In the subsequent sections, we will
propose a generalization of (a) and (b) to the case of shifted
fractional Brownian motions -- see Proposition \ref{cor_thm2}.

In order to prove Theorem \ref{thm1}, we state an easy but quite
useful lemma:
\begin{lemma}\label{derivee=cov}
Let $Z=(Z_t)_{t\in\oT}$ be a centered Gaussian process, and let
$\mathcal{H}_{1}(Z)$ be its first Wiener chaos. Fix
$Y\in\mathcal{H}_1(Z)$ and $t\in(0,T)$. Then, the following
assertions are equivalent:
\begin{enumerate}\label{equiv}
\item[(a)] \mbox{$Y$ ${\rm a.s.}$-differentiates $Z$ at $t$.}
\item[(b)] \mbox{$Y$ ${\rm L^2}$-differentiates $Z$ at $t$.}
\item[(c)] \mbox{$\frac{d}{ds}\,\cov(Z_s,Y)|_{s=t}$ exists and is finite.}
\end{enumerate}
If either (a), (b) or (c) are verified and $P(Y=0)<1$, one has
moreover that
\begin{equation}\label{cal-eq}
D^{Y}Z_t=\frac{Y}{{\rm Var}(Y)}.\frac{d}{ds}\,\cov(Z_s,Y)|_{s=t}.
\end{equation}
In particular, for every $s,t\in(0,T)$, we have: $Z_s$ ${\rm L}^2\star{\rm a.s.}$-differentiates $Z$
at $t$ if, and only if, $u\mapsto
\cov(Z_s,Z_u)$ is differentiable at $u=t$.

On the other hand, suppose that $Y\in\mathcal{H}_1(Z)$ is such that:
(i) $P(Y=0)<1$, and (ii) $Y$ does not ${\rm L}^2\star{\rm
a.s.}$-differentiate $Z$ at $t\in(0,T)$. Then, for every
$\mathscr{H}\subset\sg\{Y\}$, either $\mathscr{H}$ does not ${\rm
L}^2\star{\rm a.s.}$-differentiate $Z$ at $t$, or $\mathscr{H}$ is
such that $\esp\left[ Y \mid \mathscr{H}\right]=0$ and
$D^{\mathscr{H}}Z_t=0$.
\end{lemma}
\begin{proof}
If $Y\in\mathcal{H}_{1}(Z)\setminus\{0\}$, we have
$$\esp\left[ \Delta _{h}Z_{t}\mid Y\right]=\frac{\cov(\Delta _{h}Z_{t},Y)}{\var(Y)}\,Y.$$
The conclusions follow. \qed
\end{proof}

We now turn to the proof of Theorem \ref{thm1}:
\begin{proof}
Since $M_n$ is an invertible matrix for any $n\in I$, the
Gram-Schmidt orthonormalization procedure can be applied to
$\{Y_i\}_{i\in I}$. For this reason we may assume, for the rest of
the proof and without loss of generality, that the family $\{
Y_{i}\}_{i\in I} $ is composed of i.i.d. random variables with
common law $\mathscr{N}(0,1)$.
\begin{enumerate}
\item The first implication is an immediate consequence of Proposition
\ref{sous-bonne}. Assume now that $N<+\infty$ and that any $Y_i$,
$i=1,\ldots,N$, ${\rm L}^2$-differentiates $Z$ at
$t$. By Lemma \ref{derivee=cov}, we have in particular that
$$\frac{d}{ds}\,\cov(Z_s,Y_i)|_{s=t}$$ exists for any
$i=1,\ldots,N$. Since
\begin{equation}\label{proj2}
\esp\left[ \Delta _{h}Z_{t}\mid \mathscr{Y}\right] =\sum_{i=1}^{N}
\cov(\Delta_h Z_t,Y_i)\,Y_i
\end{equation}
we deduce that $\mathscr{Y}$ ${\rm L}^2$-differentiates $Z$ at $t$.

\item By definition, if $\mathscr{Y}$ really does not ${\rm L}^2\star{\rm a.s.}$-differentiate $Z$ at
$t$, then any finite linear combination of the $Y_i$'s either
${\rm L}^2\star {\rm a.s.}$-degenerates, or does not ${\rm
L}^2\star{\rm a.s.}$-differentiate $Z$ at $t$.

Conversely, assume that any finite linear combination of the
$Y_i$'s either ${\rm L}^2\star {\rm a.s.}$-degenerates or does not
${\rm L}^2\star{\rm a.s.}$-differentiate $Z$ at $t$. Let $\scr{G}
\subset\mathscr{Y}$. By the projection principle, we can write:
\begin{equation}\label{superproj}
{\rm E}\left[ \Delta _{h}Z_{t}\mid \scr{G}\right] =\sum_{i\in
I}\cov\left( \Delta _{h}Z_{t},Y_{i}\right) {\rm E}\left[ Y_{i}\mid
\scr{G}\right].
\end{equation}

Let us assume that $\scr{G}$ ${\rm L}^2\star{\rm
a.s.}$-differentiates $Z$ at $t$. By (\ref{superproj}) this implies
in particular that, for almost all {\it fixed} $\omega_0\in\Omega$,
$$
\esp\left[\Delta_{h}Z_t\mid\scr{G}\right](\omega_0)
=\cov\left(\Delta_{h}Z_t,\sum_{i=1}^N a_i(\omega_0)\,Y_i\right),
$$
converges as $h\rightarrow 0$, where
$a_i(\omega_0)=\esp\left[Y_i\mid\scr{G}\right](\omega_0)$. Due to
Lemma \ref{derivee=cov}, we deduce that $ X^{(\omega_0)}\triangleq
\sum_{i=1}^N a_i(\omega_0)\,Y_i$ ${\rm L}^2\star{\rm
a.s.}$-differentiates $Z$ at $t$ for almost all $\omega_0\in\Omega$.
By hypothesis, we deduce
 that $X^{(\omega_0)}$ ${\rm L}^2\star{\rm a.s.}$-degenerates $Z$ at $t$
for almost all $\omega_0\in\Omega$. But, by Lemma
\ref{derivee=cov}, the stochastic derivative
$D^{X^{(\omega_0)}}Z_t$ necessarily writes
$c(\omega_0)X^{(\omega_0)}$ with $c(\omega_0)\in\R$. Since
$X^{(\omega_0)}$ is centered and ${\rm
Var}(D^{X^{(\omega_0)}}Z_t)=0$, we deduce that
$D^{X^{(\omega_0)}}Z_t=0$. Thus
$$\lim_{h\rightarrow 0}\cov(\Delta_{h}Z_t,X^{(\omega_0)}) =\lim_{h\rightarrow
0}\esp\left[\Delta_{h}Z_t\mid\scr{G}\right](\omega_0)=0$$ for
almost all $\omega_0\in\Omega$. Thus $\scr{G}$ ${\rm
a.s.}$-degenerates $Z$ at $t$. Since $\scr{G}$ also ${\rm
L}^2$-differentiates $Z$ at $t$, we conclude that $\scr{G}$ ${\rm
L}^2\star{\rm a.s.}$-degenerates $Z$ at $t$. The proof that
$\mathscr{Y}$ really does not ${\rm L}^2\star{\rm
a.s.}$-differentiate $Z$ at $t$ is complete.

\item Again by definition, if $\mathscr{Y}$ really does not ${\rm L}^2\star{\rm a.s.}$-differentiate $Z$ at
$t$, then any finite linear combination of the $Y_i$'s either
${\rm L}^2\star {\rm a.s.}$-degenerates, or does not ${\rm
L}^2\star{\rm a.s.}$-differentiate $Z$ at $t$. We shall now assume
that every finite linear combination of the $Y_i$'s either ${\rm
L}^2\star {\rm a.s.}$-degenerates or does not ${\rm L}^2\star{\rm
a.s.}$-differentiate $Z$ at $t$. Let $(J_m)_{m\in\N}$ be the
increasing sequence given by $J_m=\{1,\ldots,m\}$, so that
$\cup_{m\in\N}
J_m=I=\N$.\\
Suppose that $\mathscr{G}\in\mathbf{R}(\scr{Y})$ and that $\scr{G}$
${\rm L}^2\star{\rm a.s.}$-differentiates $Z$ at $t$. By Proposition
\ref{sous-bonne}, $\scr{G}^i\triangleq \scr{G}\cap \sg\{Y_i\}$ ${\rm
L}^2$-differentiates $Z$ at $t$, for any $i\in\N$. But
$$
{\rm E}[\Delta_h Z_t|\scr{G}^i]={\rm
Cov}(\Delta_hZ_t,Y_i)E[Y_i|\scr{G}^i].
$$
So, for any $i\in\N$:
\begin{equation}\label{ex-or-zero}
\mbox{either }\,\lim_{h\rightarrow 0}{\rm Cov}(\Delta_h
Z_t,Y_i)\,\mbox{ exists,}\quad\mbox{or }\,{\rm
E}[Y_i|\scr{G}^i]=0.
\end{equation}
Set $\scr{G}_{m}\triangleq \scr{G}\cap \sg(Y_j, j\in J_m)$, and
observe that, if $\scr{G}\in\mathbf{R}(\scr{Y})$, then
$${\rm E}[Y_i|\scr{G}^i]={\rm E}[Y_i|\scr{G}_m]$$ for every $i=1,...,m$. We have
\begin{equation}\label{(9)}
{\rm E}\left[ \Delta _{h}Z_{t}\mid \scr{G}_m\right] =\sum_{i\in
J_m}\cov\left( \Delta _{h}Z_{t},Y_{i}\right) {\rm E}\left[
Y_{i}\mid \scr{G}^i\right].
\end{equation}
By (\ref{ex-or-zero}), and since $J_m$ is finite, we deduce that
$\scr{G}_m$ ${\rm L}^2\star{\rm a.s.}$-differentiates $Z$ at $t$.
By the same proof as in step (a) for $\scr{G}_m$ instead of
$\scr{G}$ and using (\ref{(9)}) instead of (\ref{superproj}), we
deduce that
$$
X_m^{(t)}\triangleq  D^{\mathscr{G}_{m}}Z_t=0.
$$
But, from Proposition \ref{sous-bonne}, we have:
$$D^{\mathscr{G}_{m}}Z_t ={\rm E}[D^{\mathscr{G}}Z_t|\mathscr{G}_{m}],\quad m\geq1.$$
Thus $\{X_m^{(t)},\,m\in\N\}$ is a (discrete) square integrable
martingale w.r.t. the filtration $\{\mathscr{G}_{m},\, m\in\N\}$.
So we conclude that
$$D^{\mathscr{G}}Z_t=\lim_{m\to\infty}X_m^{(t)}=0 \quad a.s.$$
In other words, $\scr{G}$ ${\rm L}^2\star{\rm a.s.}$-degenerates
$Z$ at $t$. Therefore, $\mathscr{Y}$ really does not ${\rm
L}^2\star{\rm a.s.}$-differentiate $Z$ at $t$. \qed

\end{enumerate}

\end{proof}

\textbf{Counterexample.} In what follows we show that, if
$N=+\infty$, the converse of the first point in the statement of
Theorem \ref{thm1} does not hold in general. Indeed, let $\left\{
\xi _{i}:i\geq 1\right\} $ be an infinite sequence of i.i.d.
centered standard Gaussian random variables. Let $\left\{
f_{i}:i\geq 1\right\} $ be a collection of deterministic functions
belonging to $L^{2}\left( \left[ 0,1\right] ,dt\right) $, such that
the following hold:

\begin{description}
\item[--] for every $i\geq 1$, $f_{i}\left( t\right) $ is differentiable in $%
t$ for every $t\in \left[ 0,1\right] $;

\item[--] there exists $A\in \left( 0,+\infty \right) $ such that, for every
$t\in \left[ 0,1\right] $, $\sum_{i=1}^{+\infty }f_{i}\left( t\right) ^{2}<A$%
.
\end{description}

Then, we may apply the It\^{o}-Nisio theorem (see \cite{IN}) to
deduce that there exists a Gaussian process $\{Z_{t}$ $:$ $t\in
\left[ 0,1\right] \}$ such that, a.s.-$\mathbb{P}$,
\[
\lim_{N\rightarrow +\infty }\sup_{t\in \left[ 0,1\right]
}\left\vert Z_{t}-\sum_{i=1}^{N}\xi_i f_{i}\left( t\right)
\right\vert =0\text{.}
\]%
Now suppose that the paths of $Z$ are a.s. not-differentiable for
every $t$. Then, by setting $\mathscr{Y}=\sigma( \xi _{i},\,i\geq
1)$, we obtain that $\mathscr{Y}$ does not
$\mathrm{L}^{2}\star \mathrm{a.s.}$
differentiate $Z$ at every $t$, although, for every $i\geq 1$ and every $%
t\in \left[ 0,1\right] $, $\xi _{i}$ $\mathrm{L}^{2}\star
\mathrm{a.s.}$ differentiates $Z$ at $t$. As an example, one can
consider the case
\[
f_{i}\left( t\right) =\int_{0}^{t}e_{i}\left( x\right) dx\text{, \ \ }i\geq 1%
\text{,}
\]%
where $\left\{ e_{i}:i\geq 1\right\} $ is any orthonormal basis of $%
L^{2}\left( \left[ 0,1\right] ,dx\right) $, so that the limit
process $Z$ is a standard Brownian motion. See also Kadota
\cite{Kad} for several related results, concerning the
differentiability of stochastic processes admitting a
Karhunen-Lo\`{e}ve type expansion. \qed

\section{\protect\normalsize Invariance properties of differentiating
$\sg$-fields and stochastic derivatives under equivalent changes of probability}

{\normalsize \label{sec:3} }

\smallskip
Let $Z$ be a Gaussian process, and let $\mathscr{G}\subset
\mathscr{F}$ be differentiating for $Z$. In this section we
establish conditions on $Z$ and $\mathscr{G}$, ensuring that
$\mathscr{G}$ is still differentiating for $Z$ after an equivalent
change of probability measure. As anticipated, this result will be
used to study the class of differentiating $\sg$-fields associated
with drifted Gaussian processes. Roughly speaking, we will show
that -- under adequate conditions -- one can study the stochastic
derivatives of a drifted Gaussian process by first eliminating the
drift through a Girsanov-type transformation. We concentrate on
$\sg$-fields generated by a single random variable. To achieve our
goals we will use several techniques from Malliavin calculus, as
for instance those developed by H. F\"{o}llmer (see \cite[Sec.
4]{follmer}) in order to compute the backward drift of a
non-Markovian Brownian diffusion.

Let $Z=(Z_t)_{t\in[0,T]}$ be a square integrable stochastic
process defined on a probability space $(\Omega,\mathscr{F},\bP)$.
We assume that, under an equivalent probability $\Q \sim \bP$, $Z$
is a centered Gaussian process (so that, in particular, $Z_t \in
{\rm L}^2(\bP)\cap{\rm L}^2(\Q)$ for every $t$). Let
$\mathcal{H}_{1}(Z,\Q)=\{Z(h),h\in\fr H\}$ be the first Wiener
chaos associated with $Z$ under $\Q$ (this means that the closure
is in ${\rm L}^2(\Q)$), canonically represented as an
\textit{isonormal Gaussian process} with respect to a separable
Hilbert space $(\fr H, \langle\cdot,\cdot\rangle_{\fr H})$. In
particular: (i) the space $\fr H$ contains the set $\cal E$ of
step functions on $[0,T]$, (ii) the covariance function of $Z$
under $\Q$ is given by $\rho_\Q(s,t)=\langle
\1_{[0,s]},\1_{[0,t]}\rangle_{\fr H}$, and (iii) the scalar
product $\langle
    \cdot,\cdot\rangle_{\fr H}$ verifies the general relation:
\begin{equation}\label{iso}
    \forall\  h,h'\in\fr H, \quad \langle
    h,h'\rangle_{\fr H}={\rm E}^{\Q}[Z(h)Z(h')]
\end{equation}
(note that, given $Z$, the properties (i)-(iii) completely
characterize the pair $(\fr H, \langle\cdot,\cdot\rangle_{\fr
H})$). We denote by $D$ the Malliavin derivative associated with
the process $Z$ under $\Q$ (the reader is referred to \cite{n} for
more details about these notions). The following result is an
extension of Theorem 22 in \cite{dn} to a general Gaussian
setting. Note that, in the following statements, we will
exclusively refer to the ${\rm L}^2$ topology.

\begin{thm}\label{girsa}
Fix $t\in(0,T)$ and select $g\in\fr H$ such that
$\langle\1_{[0,t]},g \rangle_\fr H\neq 0$. We write $\eta$ to
indicate the Radon-Nikodym derivative of $\Q$ with respect to
$\bP$ (that is, $d\Q=\eta\,d\bP$), and we assume that $\eta$ has
the form $\eta={\rm exp}(-\zeta)$, for some random variable
$\zeta$ for which $D\zeta$ exists. Suppose that
\begin{equation}\label{conv-girsa}
\mu_t\triangleq\lim_{h\rightarrow 0} h^{-1} \langle{\bf
1}_{[t,t+h]},D\zeta\rangle_\fr H \quad\mbox{exists in the }{\rm
L}^2 \mbox{ topology}.
\end{equation}
Then, $Z(g)$ ${\rm L}^2$-differentiates $Z$ at $t$ under $\bP$ if,
and only if, $Z(g)$ ${\rm L}^2$-differentiates $Z$ at $t$ under
$\Q$, that is, if, and only if,
\begin{equation}\label{mu}
\frac{d}{du}\,\langle{g},{\bf 1}_{[0,u]}\rangle_{\fr H}|_{u=t} =
\frac{d}{du}\,{\rm Cov}^{\Q}(Z(g),Z_u)|_{u=t} \quad\mbox{exists.}
\end{equation}
\noindent Moreover,

\begin{enumerate}
\item If $Z(g)$ ${\rm L}^2$-differentiates $Z$ at $t$ under $\Q$, then
\begin{equation}\label{conv-girsa2}
D^{Z(g)}_{\bP}Z_t=\frac{|g|^2_{\fr H}\ {\rm
E}^\bP[Z_t-\langle\1_{[0,t]},D\zeta\rangle_\fr H|Z(g)]} {Z(g)\
\langle g,\1_{[0,t]}\rangle_\fr H}\ D^{Z(g)}_{\Q}Z_t + {\rm
E}^\bP[\mu_t|Z(g)].
\end{equation}
\item If $Z(g)$ does not ${\rm L}^2$-differentiate $Z$ at $t$ under $\Q$, then
$\mathscr{H}\subset\sigma\{Z(g)\}$ differentiates $Z$ at $t$ with
respect to $\bP$ if, and only if,
$${\rm E}^\bP[Z_t-\langle\1_{[0,t]},D\zeta\rangle_\fr H|\mathscr{H}]=0.$$
In this case, $ D^\mathscr{H}_{\bP}Z_t={\rm
E}^\bP[\mu_t|\mathscr{H}]. $
\end{enumerate}
\end{thm}

\textbf{Remark.} Since $Z$ is Gaussian under $\Q$, Corollary
\ref{cor1} implies that $Z(g)$ is not differentiating for $Z$ at $t$
w.r.t. $\Q$ if, and only if, $Z(g)$ is \textsl{really} not
differentiating w.r.t. $\Q$. Point 2 in Theorem \ref{girsa} shows
that this double implication does not hold, in general, under the
equivalent probability $\bP$. Indeed, even if $Z(g)$ does not
differentiate $Z$ under $\bP$ (and therefore under $\Q$), one may
have that there exists a differentiating $\scr{H}\subset\sg\{Z(g)\}$
such that $D_{\bP}^\scr{H}Z_t$ is non-deterministic. Observe,
however, that $D_{\bP}^\scr{H}Z_t$ is forced to have the particular
form $D_{\bP}^\scr{H}Z_t= {\rm E}^\bP[\mu_t|\mathscr{H}]$.

\begin{proof}
Let $\xi\in {\rm L}^2(\bP)\cap{\rm L}^2(\Q)$ and
$A\in\mathscr{G}\subset \mathscr{F}$. The relation $$\int_A \xi
d\bP=\int_A {\rm E}^{\bP}[\xi|\scr G]d\bP$$ implies

\begin{eqnarray*}
\int_A {\rm E}^{\Q}[\xi\eta^{-1}|\scr G]d\Q &=& \int_A \xi
\eta^{-1} d\Q=\int_A
{\rm E}^{\bP}[\xi|\mathscr{G}]\eta^{-1} d\Q \\
 &=& \int_A
{\rm E}^{\bP}[\xi|\mathscr{G}]{\rm
E}^{\Q}\left[\eta^{-1}|\mathscr{G}\right]d\Q.
\end{eqnarray*}
\noindent Thus:
\begin{equation}\label{change_measure_cond_exp}
{\rm E}^{\bP}[\xi|\scr G]=\frac{{\rm
E}^{\Q}\left[\xi\eta^{-1}|\mathscr{G}\right]}{{\rm
E}^{\Q}\left[\eta^{-1}|\mathscr{G}\right]},
\end{equation}
from which we deduce that the study of ${\rm
E}^{\bP}[\Delta_hZ_t|Z(g)]$ can be reduced to that of ${\rm
E}^{\Q}[\eta^{-1}\Delta_hZ_t|Z(g)]$. Let $\phi\in
\cal{C}^1_b(\R)$. We have

\begin{eqnarray*}
  {\rm E}^{\Q}[(Z_{t+h}-Z_t)\eta^{-1} \phi(Z(g))] &=& {\rm E}^\Q[\langle\1_{[t,t+h]},D(\eta^{-1} \phi(Z(g)))\rangle_{\fr H}] \\
    &=& {\rm E}^\Q[\phi(Z(g))\eta^{-1}\langle\1_{[t,t+h]},D\zeta\rangle_{\fr H}]\\
    & & \quad +\langle\1_{[t,t+h]},g\rangle_{\fr H}
    {\rm E}^\Q[\eta^{-1}\phi'(Z(g))].
\end{eqnarray*}

\noindent By using an analogous decomposition for ${\rm
E}^\Q[Z_t\eta^{-1}\phi(Z(g))]$, we can also write:
\begin{equation}\label{}
    {\rm E}^\Q[\eta^{-1}\phi'(Z(g))]=\frac{{\rm E}^\Q[(Z_t-\langle\1_{[0,t]},D\zeta\rangle_{\fr H}
)
    \eta^{-1}\phi(Z(g))]}{\langle\1_{[0,t]},g\rangle_{\fr H}}.
\end{equation}
\noindent Therefore, ${\rm E}^\Q[\eta^{-1}\Delta_hZ_t|Z(g)]$ is
equal to
$$
{\rm E}^\Q[(Z_t-\langle\1_{[0,t]},D\zeta\rangle_{\fr H} )
    \eta^{-1}|Z(g)]
\frac{ \left\langle\1_{[t,t+h]},g\right\rangle_{\fr H} } {
h\langle\1_{[0,t]},g\rangle_{\fr H} } +
h^{-1}\langle\1_{[t,t+h]},{\rm
E}^\Q[\eta^{-1}D\zeta|Z(g)]\rangle_{\fr H},
$$
whereas ${\rm E}^\bP[\Delta_hZ_t|Z(g)]$ equals the following
expression:
\begin{equation}\label{boooh}
{\rm E}^\bP[Z_t-\langle {\bf 1}_{[0,t]},D\zeta\rangle_\fr H |Z(g)]
\frac{ \left\langle\1_{[t,t+h]},g\right\rangle_{\fr H} } {
h\langle\1_{[0,t]},g\rangle_{\fr H} } +
h^{-1}\left\langle\1_{[t,t+h]},{\rm
E}^\bP[D\zeta|Z(g)]\right\rangle_{\fr H}.
\end{equation}

\noindent Now, by assumption (\ref{conv-girsa}) and thanks to
Proposition \ref{sous-bonne}, we have that
$$
\lim_{h\rightarrow 0}h^{-1}\left\langle\1_{[t,t+h]},{\rm
E}^\bP[D\zeta|Z(g)]\right\rangle_{\fr H} ={\rm
E}^\bP[\mu_t|Z(g)]\mbox{ in the }{\rm L}^2 \mbox{  topology}.
$$
\noindent Note moreover that $\bP({\rm
E}^\bP[Z_t-\langle\1_{[0,t]},D\zeta\rangle_{\fr H}|Z(g)]=0)<1$.
Indeed, if it was not the case, one would have ($\delta$ stands
for the Skorohod integral)
\begin{eqnarray*}
0 &=& {\rm
E}^\Q[(Z_t\eta^{-1}-\langle\1_{[0,t]},D\eta^{-1}\rangle_{\fr
H})Z(g)]
={\rm E}^\Q[\delta(\1_{[0,t]}\eta^{-1})Z(g)]\\
 &=&{\rm E}^\Q[\eta^{-1}]\langle\1_{[0,t]},g\rangle_\fr H=\langle\1_{[0,t]},g\rangle_\fr H\neq
 0\,
\end{eqnarray*}
which is clearly a contradiction. As a consequence, we deduce from
(\ref{boooh}) that $Z(g)$ ${\rm L}^2$-differentiates $Z$ at $t$
under $\bP$ if, and only if, $\frac{d}{du}\,\langle{g},{\bf
1}_{[0,u]}\rangle_{\fr H}|_{u=t}$ exists. By Lemma
\ref{derivee=cov}, this last condition is equivalent to $Z(g)$
being ${\rm L}^2$-differentiating for $Z$ at $t$ under $\Q$. We
can therefore deduce (\ref{conv-girsa2}) from (\ref{boooh}) and
(\ref{cal-eq}).

\noindent If $\mathscr{H}\subset\sg\{Z(g)\}$, the projection
principle and (\ref{boooh}) yield that ${\rm
E}^\bP[\Delta_hZ_t|\mathscr{H}]$ equals
$$
{\rm E}^\bP[Z_t-\langle {\bf 1}_{[0,t]},D\zeta\rangle_\fr H
|\mathscr{H}] \frac{ \left\langle\1_{[t,t+h]},g\right\rangle_{\fr
H} } { h\langle\1_{[0,t]},g\rangle_{\fr H} } +
h^{-1}\left\langle\1_{[t,t+h]},{\rm
E}^\bP[D\zeta|\mathscr{H}]\right\rangle_{\fr H}.
$$
\noindent When $\frac{d}{du}\,\langle{\bf 1}_{[0,s]},{\bf
1}_{[0,u]}\rangle_{\fr H}|_{u=t}$ does not exist, we deduce that
$\mathscr{H}$ differentiates $Z$ at $t$ if, and only if, ${\rm
E}^\bP[Z_t-\langle {\bf 1}_{[0,t]},D\zeta\rangle_\fr H
|\mathscr{H}]=0$. If this condition is verified, we then have $
D^\mathscr{H}_{\bP}Z_t={\rm E}^\bP[\mu_t|\mathscr{H}], $ again by
Proposition \ref{sous-bonne}. \qed
\end{proof}

As an application of Theorem \ref{girsa}, we shall consider the
case where the isonormal process $Z$ in (\ref{iso}) is generated
by a fractional Brownian motion of Hurst index
$H\in(0,1/2)\cup(1/2,1)$ (see also \cite[Theorem 22]{dn}, for
related results concerning the case $H\in(1/2,1)$).

We briefly recall some basic facts about stochastic calculus with
respect to a fractional Brownian motion. We refer the reader to
\cite{n-cours} for any unexplained notion or result. Let
$B=(B_t)_{t\in[0,T]}$ be a fractional Brownian motion with Hurst
parameter $H\in (0,1)$, and assume that $B$ is defined on a
probability space $(\Omega, \scr F,\bP)$. This means that $B$ is a
centered Gaussian process with covariance function ${{\rm E}}
(B_sB_t)=R_H(s,t)$ given by
\begin{align}
\label{cov} R_H(s,t)=\frac{1}2\left (
t^{2H}+s^{2H}-|t-s|^{2H}\right).
\end{align}
We denote by $\mathcal E$ the set of all $\R-$valued step
functions on $\oT$. Let $\fr H$ be the Hilbert space defined as
the closure of $\mathcal E$ with respect to the scalar product
$$
\left \langle\1_{[0,t]}, g\right \rangle_{\fr H}=R_H(t,s),
$$
and denote by $|\cdot|_{\fr H}$ the associate norm. The mapping
$\displaystyle \1_{[0,t]} \mapsto B_{t}$ can be extended to an
isometry between $\fr H$ and the Gaussian space $\cal H_1(B)$
associated with $B$. We denote this isometry by $\ffi\mapsto
B(\ffi)$. Recall that the covariance kernel $R_H(t,s)$ introduced
in (\ref{cov}) can be written as
$$ R_H(t,s)= \int_{0}^{s\wedge t} K_H(s,u)K_H(t,u)du,$$
where $K_H(t,s)$ is the square integrable kernel defined, for $s<t$,
by
\begin{equation}\label{du}
K_{H}(t,s)=\Gamma(H+\frac{1}{2})^{-1}(t-s)^{H-\frac{1}{2}}F\big(H-\frac{1}{2},\frac{1}{2}-H,H+\frac{1}{2},1-\frac{t}{s}\big),
\end{equation}
where $F(a,b,c,z)$ is the classical Gauss hypergeometric function.
By convention, we set $K_H(t,s)=0$ if $s\ge t$. We define the
operator $\mathcal K_H$ on ${\rm L}^2([0,T])$ as
$$
(\cal K_H h)(t)=\int_0^tK_H(t,s)h(s)ds. $$ Let $\mathcal
K_{H}^{\ast}:\mathcal E\to {\rm L}^2([0,T])$
be the linear operator defined as: %
$$
\mathcal K_{H}^{\ast}\left (  \1_{[0,t]}\right ) = K_H(t,\cdot).
$$
The following equality holds for any $\phi,\psi \in\mathcal E$
$$
\langle \phi,\psi\rangle_{\fr H} = \langle \mathcal K_{H}^{\ast}
\phi,\mathcal K_{H}^{\ast} \psi\rangle_{{\rm L}^{2}([0,T])}={{\rm
E}}\left ( B(\phi)B(\psi)\right ),
$$
implying that $\mathcal K_{H}^{\ast}$ is indeed an isometry
between the Hilbert spaces $\fr H$ and a closed subspace of ${\rm
L}^{2}([0,T])$. Now consider the process $W = (W_t)_{t\in [0,T]} $
defined as
$$
W_t  = B\big ((\mathcal K_H^\ast)^{-1}(\1_{[0,t]})\big ),
$$
and observe that $W$ is a standard Wiener process, and also that
the process $B$ has an integral representation of the type
$$
B_t =\int_{0}^t K_H(t,s) dW_s,
$$
so that, for any $\phi\in\fr H$,
$$
B(\phi)  = W\left ( \mathcal K_H^\ast \phi \right ).
$$

We will also need the fact that the operator ${\cal K}_H$ can be
expressed in terms of fractional integrals as follows:
\begin{eqnarray}
\label{kh1}({\cal K}_H h)(s)&=&I_{0+}^{2H} s^{\frac{1}{2}-H} I_{0+}^{\frac{1}{2}-H}s^{H-\frac{1}{2}}h(s),\quad\mbox{if $H< 1/2$},\\
\label{kh2}({\cal K}_H h)(s)&=&I_{0+}^{1} s^{H-\frac{1}{2}}
I_{0+}^{H-\frac{1}{2}}s^{\frac{1}{2}-H}h(s),\quad\mbox{if $H> 1/2$},
\end{eqnarray}
for every $h\in{\rm L}^2([0,T])$. Here, $I_{0+}^\alpha f$ denotes
the left fractional Riemann-Liouville integral of order $\alpha$
of $f$, which is defined by
$$
I_{0+}^\alpha f(x)=\frac{1}{\Gamma(\alpha)}\int_0^x
(x-y)^{\alpha-1}f(y)dy.
$$

Let $\Upsilon_H$ be the set of the so-called shifted fBm
$Z=(Z_t)_{t\in [0,T]}$ defined by
\begin{equation}\label{Z}
Z_t=x_0+ B_t + \int_0^t b_s ds,\quad t\in [0,T],
\end{equation}
where $b$ runs over the set of adapted processes (w.r.t. the natural filtration of $B$) having
integrable trajectories.

\smallskip
We also need to introduce a technical assumption. Define
\begin{equation}\label{jj}
a_r  =  \left(\cal K_H^{-1}\int_0^{\cdot}b_s ds\right)(r);
\end{equation}
in what follows we shall always assume that
\begin{itemize}
\item[(H1)] $a$ is bounded a.s.,
\item[(H2)] $\Phi$ defined by $\Phi(s)=\int_0^T D_s a_r \delta B_r$
exists and belongs in ${\rm L}^2([0,T])$ a.s..
\end{itemize}

First, let us consider the case $H>1/2$. We suppose moreover that the
trajectories of $b$ are a.s. H\"older continuous of order
$H-1/2+\e$, for some $\e>0$. Then, the fractional version of the
Girsanov theorem (see \cite[Theorem 2]{no}) applies, yielding that
$Z$ is a fractional Brownian motion of Hurst parameter $H$ under
the new probability $\Q$ defined by $d\Q=\eta d\bP$, where
\begin{equation}\label{etaH}
\eta={\rm exp}\left(-\int_0^T \big({\cal K}_H^{-1}\int_0^\cdot
b_rdr\big)(s)dW_s -\frac{1}{2}\int_0^T \big({\cal
K}_H^{-1}\int_0^\cdot b_rdr\big)^2(s)ds\right).
\end{equation}
We can now state the following extension of Theorem 22 in
\cite{dn}:
\begin{cor}\label{cor_thm2}
Let $Z\in\Upsilon_H$ with $H>1/2$ and $s,t\in(0,T)$. Then $Z_s$
${\rm L}^2$-differentiates $Z$ at $t$.
\end{cor}
\begin{proof}
The proof of this result relies on Theorem \ref{girsa}. Note also
that parts of the arguments rehearsed below are only sketched,
since they are analogous to those involved in the proof of
\cite[Theorem 22]{dn}. Let us consider
$$
\zeta=\int_0^Ta_sdW_s+\frac{1}{2}\int_0^Ta_s^2ds,
$$
where $a$ is defined according to (\ref{jj}). We shall show that
(\ref{conv-girsa}) holds. We can compute (see the proof of
\cite[Theorem 22]{dn})
\begin{equation}\label{dn-calcul}
\langle{\bf 1}_{[t,t+h]},D\zeta\rangle_\fr H =\int_t^{t+h}b_rdr +
({\cal K}_H\Phi)(t+h)-({\cal K}_H\Phi)(t),
\end{equation}
where $\Phi(s)=\int_0^T D_s a_r \delta B_r$, see (H2). Since, in
the case where $H>1/2$, ${\cal K}_H\Phi$ is differentiable at $t$
(see for instance (\ref{kh2})) we deduce that (\ref{conv-girsa})
holds. Moreover, one can easily prove that (\ref{mu}) also holds,
so that the proof is concluded. \qed
\end{proof}

Now we consider the case $H<1/2$. We assume moreover that
$\int_0^T b_r^2dr<+\infty$ a.s.. Then, the fractional version of
the Girsanov theorem (see \cite[Theorem 2]{no}) holds again,
implying that $Z$ is a fractional Brownian motion of Hurst
parameter $H$ under the new probability $\Q$ defined by $d\Q=\eta
d\bP$, with $\eta$ given by (\ref{etaH}). Note that, when $H<1/2$,
we cannot apply Theorem \ref{girsa}, since (\ref{conv-girsa}) does
not hold in general. The reason is that ${\cal K}_H\Phi$ is no
more differentiable at $t$, see (\ref{kh1}). In order to make
$h^{-1}{\rm E}[Z_{t+h}-Z_t|Z(g)]$ converge, we have to replace
$h^{-1}$ with $h^{-2H}$ (we only consider the case where $h>0$).
This fact is made precise by the following result.
\begin{pro} \label{prop<1/2}
Let $Z\in\Upsilon_H$ with $H<1/2$ and $s,t\in(0,T)$. Then,
$$
\lim_{h\downarrow 0} h^{-2H}{\rm E}[Z_{t+h}-Z_t|Z_s]\mbox{ exists
in the ${\rm L}^2$-topology}.
$$
\end{pro}
\begin{proof}
We go back to the proof of Corollary \ref{cor_thm2}, with special
attention to relation (\ref{dn-calcul}). By setting
$\phi(s)=s^{\frac{1}{2}-H}I_{0+}^{\frac{1}{2}-H}s^{H-\frac{1}{2}}\Phi(s)$,
we have
\begin{eqnarray*}
{\cal K}_H\Phi(t+h)-{\cal K}_H\Phi(t)
&=&
I_{0+}^{2H}\phi(t+h) - I_{0+}^{2H}\phi(t)\\
&=&\frac{1}{\Gamma(2H)}\int_0^t \big( (t+h-y)^{2H-1}-(t-y)^{2H-1}\big)\phi(y)dy\\
&&\quad + \frac{1}{\Gamma(2H)}\int_t^{t+h} (t+h-y)^{2H-1}\phi(y)dy\\
&=&\frac{1}{\Gamma(2H)}\int_0^t \big( (y+h)^{2H-1}-y^{2H-1}\big)\phi(t-y)dy\\
&&\quad + \frac{1}{\Gamma(2H)}\int_0^{h} y^{2H-1}\phi(t+h-y)dy\\
&=&\frac{h^{2H}}{\Gamma(2H)}\int_0^{t/h} \big( (y+1)^{2H-1}-y^{2H-1}\big)\phi(t-hy)dy\\
&&\quad +\frac{h^{2H}}{\Gamma(2H)}\int_0^{1} y^{2H-1}\phi(t+h-hy)dy.
\end{eqnarray*}
We deduce that
$$
h^{-2H}\big({\cal K}_H\Phi(t+h)-{\cal
K}_H\Phi(t)\big)\longrightarrow c_H\,\phi(t), \mbox{ as
$h\rightarrow 0$},
$$
where
$$
c_H=\frac{1}{\Gamma(2H)}\int_0^{+\infty}\big((y+1)^{2H-1}-y^{2H-1}\big)dy+\frac{1}{\Gamma(2H)}\int_0^1
y^{2H-1}dy<+\infty.
$$
Thus, by using the notations adopted in (the proof of) Theorem
\ref{girsa}, one deduces an analogue of (\ref{conv-girsa}),
obtained by replacing $h^{-1}$ with $h^{-2H}$, that is:
$$
\tilde{\mu}_t\triangleq\lim_{h\rightarrow 0} h^{-2H} \langle{\bf
1}_{[t,t+h]},D\zeta\rangle_\fr H \quad\mbox{exists in the }{\rm
L}^2 \mbox{ topology}.
$$
Moreover, it is easily shown that  $\lim_{h\rightarrow
0}h^{-2H}\langle{\bf 1}_{[t,t+h]},{\bf 1}_{[0,s]}\rangle_{\fr H}$
exists. By using (\ref{boooh}), we obtain the desired conclusion.
\qed
\end{proof}

\section{Differentiating collections of $\sigma$-fields and the
associated differentiated process}

In this section, we work on a complete probability space
$(\Omega,\scr F, \bP)$, and we denote by $\scr B_{(0,T)}$ the
Borel $\sg$-field of $(0,T)$. In the previous sections, we have
studied the properties of those $\sg$-field that are
differentiating for some processes at a \textsl{fixed} time $t$.
We will now concentrate on collections of differentiating
$\sg$-fields indexed by the whole interval $(0,T)$.

\begin{defi}
{\normalsize We say that a collection $(\mathscr{A}^t)_{t\in(0,T)}$
of $\sigma$-fields
$\tau$-differentiates $Z$ if, for any $%
t\in(0,T)$, $\mathscr{A}^t$ $\tau$-differentiates $Z$ at $t $. }
\end{defi}

A differentiating collection of $\sigma$-fields need not be a
filtration (see {\it e.g.} section 5 in \cite{dn}). Nevertheless,
we can associate to each $\tau$-differentiating collection $\fr
A=(\scr{A}^t)_{t\in(0,T)}$ for $Z$ a filtration $\cal
A=(\scr{A}_t)_{t\in(0,T)}$, obtained by setting:
$$\scr{A}_t=\bigvee_{0< s\leq t}\scr A^t,\qquad t\in(0,T).$$

The collection of r.v. $(D^{\scr A^t}Z_t)_{t\in (0,T)}$ is a $\cal
A$-adapted process \cite[Definition 27.1]{rogers_williams}, in the
sense that for all $t\in (0,T)$, $D^{\scr A^t}Z_t$ is $\scr
A_t$-measurable. We call it the \textit{differentiated process} of
$Z$ w.r.t.
$\fr A$, and we denote it by $D^{\fr A}Z$.\\

In order to use such a process in stochastic analysis, one should
know whether it admits a measurable version, that is, whether
there exists a process $Y$ which is $\scr B_{(0,T)}\otimes \scr
F$-measurable and such that for all $t$, $Y_t=D^{\scr A^t}Z_t$
a.s.. Our aim in this section is to obtain a sufficient condition
for the existence of a measurable version. To this end, we
introduce the following

\begin{defi}\label{}
Let $\fr A=(\scr A^t)_t$ be a collection of $\sg$-fields and $Z$ be
a measurable stochastic process. We say that $\fr A$ is {\rm
regular} for $Z$ if for all $n\in \N$, $i\in\{1,\cdots,n\}$,
$t_i\in[0,T]$, $\phi_i\in C^{\infty}_0(\R^d)$, the process
$$t\mapsto {\rm E}[\phi_1(Z_{t_1})\cdots\phi_n(Z_{t_n})|\scr A^t]$$
has a measurable version.
\end{defi}

If $\fr A$ is a filtration, then $\fr A$ is regular for any process.
For Gaussian processes and most of drifted Gaussian processes $X$,
the collection $\fr A=(\sg\{X_t\})_t$ is a regular
collection for $X$.\\

The next result shows that, under the {\it regularity} condition
defined above, a measurable version of the differentiated process
exists. This follows from one of Doob's most celebrated theorems
(see \textit{e.g.} \cite[Theorem 30 p.158]{dm}).

\begin{thm}\label{version_mes}
Let $X$ be a $\scr B_{(0,T)}\otimes \scr F$-measurable stochastic
process defined on a complete probability space $(\Omega,\scr
F,\bP)$, and assume that $\scr F=\sigma\{X\}$. Let $\fr A$ be a
regular $L^1$-differentiating collection for $X$. Then, there exists
a measurable version of the differentiated process $D^{\fr A}X$.
This version is also adapted to the filtration generated by $\fr A$.
\end{thm}

\begin{proof}\label{}
Fix $\e >0$, and let $(h_k)$ be a sequence converging to $0$ and
$Z^k$ be the process defined by
$$Z^k_t=\frac{X_{t+h_k}-X_t}{h_k}.$$ Since $X$ is measurable, so
is the process $Z$. Then, by \cite[Theorem 30 p.158]{dm}, there
exist elementary processes $U^{n_k}_t$ such that, for all
$t\in(0,T)$ and every $k$, ${\rm E}|Z^k_t-U^{n_k}_t|<\e/2$. These
elementary processes have the form:
$$U^{n_k}_t=\sum_i
\1_{A^{n_k}_i}(t)H^{n_k}_i,$$ where $(A^{n_k}_i)_i$ is a finite
partition of $(0,T)$ and $H^{n_k}_i$ are $\scr F$-measurable
random variables. We have
$${\rm E}[U^{n_k}_t|\scr A^t]=\sum_i
\1_{A^{n_k}_i}(t){\rm E}[H^{n_k}_i|\scr A^t].$$ Since cylindrical
functionals of $X$ are dense in $L^1(\Omega, \scr F)$, we deduce
from the regularity condition that the processes $t\mapsto {\rm
E}[H^{n_k}_i|\scr A^t]$ admits a $\scr B_{(0,T)}\otimes \scr
F$-measurable modification and also, by linearity, the same
conclusion holds for the process $t\mapsto {\rm E}[U^{n_k}_t|\scr
A^t]$. Moreover,
$${\rm E}\left|{\rm E}[Z^k_t|\scr A^t]-{\rm E}[U^{n_k}_t|\scr A^t]\right|<\e/2.$$
Since $\fr A$ is a $L^1$-differentiating collection for $X$, we
deduce that there exists $k$ such that $${\rm E}\left|D^{\scr
A^t}X_t-{\rm E}[Z^k_t|\scr A^t]\right|<\e/2,$$ and therefore ${\rm
E}\left|D^{\scr
A^t}X_t-E[U^{n_k}_t|\scr A^t]\right|<\e$ for every $t$.\\
We now deduce that the map $t\mapsto [D^{\scr A^t}X_t]$ is
measurable, where $[\cdot]$ denotes the class of a process in
$L^1(\Omega)$ reduced by null sets. Indeed, it is the limit in the
Banach space $L^1(\Omega)$ (when $k$ goes to infinity) of the
measurable map $t\mapsto {\rm E}[U^{n_k}_t|\scr A^t]$. Since
$L^1(\Omega)$ is separable, we again deduce from \cite[Theorem 30
p.158]{dm} that $D^{\fr A}X$ admits a measurable modification.
\qed
\end{proof}

\section{\protect\normalsize Embedded differential equations}

{\normalsize \label{emb} }

The last section of the paper is devoted to the outline of a
general framework for \textit{stochastic embedding problems}
(introduced in \cite{cd}) related to ordinary differential
equations. As we will see, this notion involves the stochastic
derivative operators that we have defined and studied in the
previous sections. Roughly speaking, the aim of a stochastic
embedding procedure is to write a "stochastic equation" which
admits both stochastic and deterministic solutions, in such a way
that the deterministic solutions also satisfy a fixed ordinary
differential equation (see \cite{cd}). It follows that the
embedded stochastic equation is a genuine extension of the
underlying ordinary differential equation to a stochastic
framework.

\subsection{\protect\normalsize General setting}

{\normalsize Let $\chi:\mathbb{R}^d\rightarrow\mathbb{R}^d$
($d\in\N^*$) be a smooth vector field. Consider the ordinary
differential equation:
\begin{equation}  \label{edo}
\frac{dx}{dt}(t)=\chi(x(t)),\quad t\in [0,T].
\end{equation}

Let $\Lambda$ be a set of measurable stochastic processes $%
X:\Omega\times [0,T]\rightarrow\mathbb{R}^d$, where $(\Omega,\scr
F,\bP)$ is a fixed probability space. In order to distinguish two
different kinds of families of $\sigma$-fields, we shall adopt the
following notation: (i) the symbol $\fr
A_0=(\mathscr{A}_0^t)_{t\in[0,T]}$ denotes a collection of
$\sg$-fields whose definition does not depend on the choice of $X$ in the class
$\Lambda$, and (ii) $\fr A=(\mathscr{A}_X^t)_{X \in \Lambda,
t\in[0,T]}$ indicates a generic family of $\sg$-fields such that,
for every $t\in[0,T]$ and every $X \in \Lambda$, $\scr
A_X^t\subset \scr P^X_T$. We introduce the following natural
assumption:
\begin{enumerate}  \label{H1}
\item[$(T)$] $\Lambda$ contains all the deterministic
differentiable functions $f:[0,T]\rightarrow\R^d$ (viewed as
deterministic stochastic processes).
\end{enumerate}

We now fix a topology $\tau$, and describe two \textit{stochastic
embedded equations} associated with (\ref{edo}).

\begin{defi}
Fix a class of stochastic processes $\Lambda$ on $(\Omega,\scr
F,\bP)$, verifying assumption (T).

{\rm (a) Given} a family $\fr
A_0=(\mathscr{A}_0^t)_{t\in[0,T]}$ of $\sg$-fields, we say that
the equation
\begin{equation}  \label{embedding1}
X\in\Lambda,\quad D^{\mathscr{A}_0^{t}}X_t=\chi(X_t)\quad\mbox{for
every } t\in [0,T],
\end{equation}
is the {\rm strong stochastic embedding} in $\Lambda$ of the ODE
(\ref{edo}) w.r.t. $\fr A_0$.

{\rm (b)} Given a family  $\fr
A=(\mathscr{A}_X^t)_{X\in\Lambda,t\in[0,T]}$ of $\sg$-fields such
that for all $X\in\Lambda$ and $t\in[0,T]$, $\scr A_X^t\subset
\scr P^X_T$, we say that the equation
\begin{equation}  \label{embedding2}
X\in\Lambda,\quad D^{\mathscr{A}_X^{t}}X_t=\chi(X_t)\quad\mbox{for
every } t\in [0,T],
\end{equation}
is the {\rm weak stochastic embedding} in $\Lambda$ of the ODE
(\ref{edo}) w.r.t. $\fr A$.

{\rm (c)} A {\rm solution} of (\ref{embedding1}) (resp.
(\ref{embedding2})) is a stochastic process $X \in \Lambda$ such
that: {\rm (c-1)} the process $D^{\mathscr{A}_0^{t}}X_t$ (resp.
$D^{\mathscr{A}_X^{t}}X_t$) admits a jointly measurable version,
and {\rm (c-2)} the equation $D^{\mathscr{A}_0^{t}}X_t=\chi(X_t)$
(resp. $D^{\mathscr{A}_X^{t}}X_t=\chi(X_t)$) is verified for every
$t \in [0,T]$.
\end{defi}

Note that a solution of (\ref{edo}) is always a solution of
(\ref{embedding1}) or (\ref{embedding2}). Observe also that if one
wants to obtain "genuinely stochastic" solutions of (\ref{edo})
({\it i.e.} non deterministic), the previous definition implicitly
imposes some restrictions on the class $\Lambda$. Namely, if $X
\in \Lambda$ is a solution of (\ref{embedding1}) (resp.
(\ref{embedding2})), then for any $t\in[0,T]$, $\mathscr{A}_0^t$
(resp. $\mathscr{A}_X^t$) is differentiating for $X$ at $t$ with
respect to the topology $\tau$ and the random variable $\chi(X_t)$
is $\mathscr{A}_0^t$-measurable (resp.
$\mathscr{A}_X^t$-measurable) for every $t$. As an example, let
$\Gamma$ be the set of all processes $X$ with the form:
\begin{equation}\label{}
    X_t=X_0+\sg B_t+\int_0^tb_rdr, \quad t\in[0,T]
\end{equation}
where $\sg \in \R$, $B$ is a fBm of Hurst index
$H\in(0,1)$, and $b$ runs over the set of adapted processes (w.r.t. the
natural filtration of $B^H$) having a.s. integrable trajectories.
Suppose that we seek for solutions with $\sg\neq 0$ of the weak
stochastic embedding of (\ref{edo}) given by
\begin{equation}\label{}
X\in\Gamma,\quad D^{\sg\{X_t\}}X_t=\chi(X_t),
\end{equation}
Then, Corollary \ref{cor_thm2} and Proposition \ref{prop<1/2}
imply that such solutions must necessarily be driven by a fBm of
Hurst index
$H>1/2$.\\

Stochastic embedded equations may be useful in the following
framework. Suppose that a physical system is described by
(\ref{edo}), and that we want to enhance this deterministic
mathematical model in order to take into account some "stochastic
phenomenon" perturbing the system. Then, the embedded equations
(\ref{embedding1}) or (\ref{embedding2}) may be the key to define
a stochastic model in a very coherent way, in the sense that every
stochastic process satisfying (\ref{embedding1}) or
(\ref{embedding2}) is also constrained by the physical laws ({\it
i.e.} the ODE (\ref{edo})) defining the original deterministic
description of the system.

\subsection{\protect\normalsize A first example}

{\normalsize Consider the set $\Lambda$ of all continuous
processes defined on the probability space
$(\Omega,\mathscr{F},\bP)$, as well as the "constant" collection
of $\sigma$-fields $(\mathscr{F}_t)_{t\in{[0{,}T]}}$ such that
$\mathscr{F}_t=\mathscr{F}$ for every $t$. Since the stochastic
derivative w.r.t. $\mathscr{F}$ coincides with the usual pathwise
derivative, the embedding problem
\begin{equation}
D^\mathscr{F}X_t=\chi(X_t),\quad t\in [0,T],
\end{equation}
has a unique strong solution for a given initial condition (deterministic
or random). Note that in this example the embedded differential
equation produces no other solution than those given by
(\ref{edo}). }

\subsection{\protect\normalsize A more interesting example}

{\normalsize Let $W$ be a Wiener process on $[0,T]$ and consider the
set $\Lambda$ of deterministic processes and of all stochastic
processes that can be expressed in terms of multiple stochastic
integrals with respect to $W$. More precisely,
denote by $\Lambda_W$ the set of processes $u\in \mathrm{L}^2(\Omega,%
\mathrm{L}^2([0,T]))$ such that, for every $t\in[0,T]$,
\begin{equation*}
u_t=\sum_{n\ge 0} J_n(f_n(\cdot,t))
\end{equation*}
where, for any $t\in[0,T]$, the $f_n(\cdot,t)$'s verify:
\begin{equation*}
\sum_{n\ge 0} \left( \left\|
f_n(\cdot,t)\right\|^2_{\mathrm{L}^2(\Delta_n[0,T])}
+\left\| \frac{\partial f_n}{\partial t}(\cdot,t)\right\|^2_{\mathrm{L}%
^2(\Delta_n[0,t])} \right)<+\infty.
\end{equation*}
Here
\begin{equation*}
\Delta_n[0,T]=\{(s_1,\ldots,s_n)\in\mathbb{R}_+^n:\,0\le s_n\le
\ldots\le s_1\le T\}
\end{equation*}
and, for $g\in \mathrm{L}^2(\Delta_n[0,T])$,
\begin{equation*}
J_n(g)=\int_{\Delta_n[0,T]}g\,dW=\int_0^T
dW_{s_1}\int_0^{s_1}dW_{s_2}\ldots\int_0^{s_{n-1}}dW_{s_{n}}g(s_1,%
\ldots,s_n).
\end{equation*}
On $\Lambda_W$, we can consider stochastic derivatives of Nelson
type ({\it i.e.} w.r.t. a fixed filtration \cite{nels}): }

\begin{lemma}
{\normalsize Fix $t\in ]0,T[$ and let $\mathscr{P}_t$ be the past before $t$%
, that is the $\sigma$-field generated by $\{W_s,\ 0\le s\le t\}$.
If $u\in \Lambda_W$ then $D^{\mathscr{P}_t}u_t$ exists and it is
given by
\begin{equation}  \label{lm-der}
D^{\mathscr{P}_t}u_t=\sum_{n\ge 0} J_n\left(\frac{\partial f_n}{%
\partial t}(\cdot,t)\mathbf{1}_{\Delta_n[0,t]}\right)\mbox{ in the ${\rm L}^2$ sense}.
\end{equation}
}
\end{lemma}
{\normalsize
\begin{proof} Obvious by projection. \qed
\end{proof}

As an example, consider the case where $\chi$ is given by
$\chi(x)=ax+b$ with $a,b\in\mathbb{R}$. In other words, we want to
solve the strong embedding
\begin{equation}  \label{ex-eq}
X\in \Lambda,\quad D^{\mathscr{P}_t}X_t=aX_t+b,\quad t\in [0,T]
\end{equation}
in the class $\Lambda_W$. It is easy to see that if $X \in
\Lambda_W$, then $X$ satisfies (\ref{ex-eq}) if, and only if, the
kernels in its chaotic expansion satisfy
\begin{equation*}
\frac{\partial f_n}{\partial t}(\cdot,t)\mathbf{1}_{\Delta_n[0,t]}(\cdot)=a%
\,f_n(\cdot,t),\quad t\in]0,T[
\end{equation*}
for any $n\in\mathbb{N}^*$ and
\begin{equation*}
f^{\prime }_0(t)=a\,f_0(t)+b,\quad t\in]0,T[.
\end{equation*}
We deduce that $X\in\Lambda_W$ solves strongly (\ref{ex-eq}) if, and only if,
there exists a sequence $(c_n)_{n\in\mathbb{N}}$ of functions from
$\Delta_n[0,T]$ to $\mathbb{R}$ such that
\begin{equation*}
f_n(\cdot,t)=c_n(\cdot)\mathrm{e}^{at}\mathbf{1}_{\Delta_n[0,t]}(\cdot),%
\quad t\in[0,T]
\end{equation*}
for every $n\in\mathbb{N}^*$, and
\begin{equation*}
f_0(t)=c_0\,\mathrm{e}^{at}-b/a,\quad t\in[0,T].
\end{equation*}\\

Several properties of embedded stochastic equations will be
investigated in a separate paper. For instance, we will be
interested in establishing conditions ensuring that the solution
of an embedded equation is Markovian. Also, we will explore
embedded stochastic equations that are obtained from ordinary
equations of order greater than one.


\end{document}